\documentclass[12pt]{article}
\usepackage{amssymb}
\usepackage{booktabs}
\def \qed {\hfill $\boxempty$}

\def \sst {\subset}
\def \ssq {\subseteq}
\def \smin {\setminus}

\def \tur {\mbox{\rm ex}}

\newtheorem{Theorem}{Theorem}
\newtheorem{lem}[Theorem]{Lemma}
\newtheorem{defi}[Theorem]{Definition}
\newtheorem{crl}[Theorem]{Corollary}
\newtheorem{prp}[Theorem]{Proposition}
\newtheorem{prm}[Theorem]{Problem}
\newtheorem{cnj}[Theorem]{Conjecture}
\newtheorem{rmk}[Theorem]{Remark}
\newtheorem{xmp}[Theorem]{Example}
\newtheorem{cl}{Claim}

\def \bcl {\begin{cl} \ }
\def \ecl {\end{cl}}

\def \bp {\begin{prp} \ }
\def \ep {\end{prp}}

\def \bc {\begin{crl} \ }
\def \ec {\end{crl}}

\def \thm {\begin{Theorem} \ }
\def \ethm {\end{Theorem}}

\def \bl {\begin{lem} \ }
\def \el {\end{lem}}

\def \bd {\begin{defi} \ \rm }
\def \ed {\end{defi}}

\def \brm {\begin{rmk} \ }
\def \erm {\end{rmk}}

\def \bxm {\begin{xmp} \ \rm }
\def \exm {\end{xmp}}

\def \bpm {\begin{prm} \ }
\def \epm {\end{prm}}

\def \bcj {\begin{cnj} \ }
\def \ecj {\end{cnj}}

\def \nmr {\begin{enumerate}}
\def \enmr {\end{enumerate}}

\def \tmz {\begin{itemize}}
\def \etmz {\end{itemize}}

\def \nin {\noindent}
\def \bsk {\bigskip}

\def \pf {\nin{\bf Proof } \ }
\def \qed {\hfill $\Box$}

\def \kp {\,:\,}

\def\cF{{\cal F}}
\def\cH{{\cal H}}
\def\cI{{\cal I}}
\def\cK{{\cal K}}
\def\cP{{\cal P}}
\def\cR{{\cal R}}
\def\frF{{\mathfrak{F}}}

\def \ext {\mbox{\rm ex}}

\begin{document}

\title{Partition-crossing
  hypergraphs~\thanks{
  ~Research supported in part by the National Research, Development and
Innovation Office -- NKFIH under the grant SNN 116095.
}}
\author{Csilla Bujt\'as~$^1$\qquad \vspace{1ex}
        Zsolt Tuza~$^{1,2}$
\\
\normalsize $^1$~Faculty of Information Technology,
 University of Pannonia\\
\normalsize H--8200 Veszpr\'em, \vspace{1ex} Egyetem u.~10,
 Hungary \\
\normalsize $^2$~Alfr\'ed R\'enyi Institute of Mathematics,
        Hungarian Academy of Sciences \\
\normalsize H--1053 Budapest,
 Re\'altanoda u.~13--15, Hungary
~ }
\date{\normalsize\sl Dedicated to the memory of our friend and colleague Csan\'ad Imreh}
\maketitle

\begin{abstract}

For a finite set $X$, we
 say that a set $H\subseteq X$ crosses a partition $\cP=(X_1,\dots,X_k)$
 of $X$ if $H$ intersects $\min (|H|,k)$ partition classes.
If $|H|\geq k$, this means that $H$ meets all classes $X_i$,
 whilst for $|H|\leq k$ the elements of the crossing set $H$ belong to mutually distinct classes.
A set system $\cH$ crosses $\cP$, if so does some $H\in \cH$.
 The minimum number of $r$-element subsets, such that every
 $k$-partition of an $n$-element set $X$ is crossed by at least one
 of them, is denoted by $f(n,k,r)$.

  The problem of determining these minimum values for $k=r$ was raised
 and studied by several authors, first  by Sterboul in 1973
  [\textit{Proc.\ Colloq.\ Math.\ Soc.\ J. Bolyai}, Vol.\ 10, Keszthely 1973,
 North-Holland/American Elsevier, 1975, pp.\ 1387--1404].
  The present authors determined asymptotically tight
  estimates on $f(n,k,k)$ for every fixed $k$ as $n\to\infty$
  [\textit{Graphs Combin.}, 25 (2009),  807--816].
 Here we consider the more general problem for two parameters $k$ and $r$,
  and establish lower and upper bounds
 for $f(n,k,r)$. For various combinations of the three values
  $n,k,r$ we obtain asymptotically tight estimates, and also point out
  close connections of the function $f(n,k,r)$ to
  Tur\'an-type extremal problems on graphs and hypergraphs,
  or to balanced incomplete block designs.

\bigskip

\noindent {\bf Keywords\kp}
 Partition, Set system, Crossing set, Tur\'an-type problem, Hypergraph,
 Upper chromatic number.

\bigskip

\nin {\bf Mathematics 2000 Subject Classification\kp}
 05C35, 
 05C15, 
05C65 

\end{abstract}

\newpage

\section{Introduction}

Let $X$ be a finite set.
By a \emph{$k$-partition}   of $X$ we mean a partition
 $\cP=(X_1,\dots,X_k)$ into \emph{precisely\/ $k$ nonempty
 classes}.
For a natural number $r\ge 2$, the family of all $r$-element
 subsets of $X$ --- also termed {\em $r$-subsets\/}, for short (similarly,
 `$r$-set' may abbreviate `$r$-element set') ---
 is denoted by $X\choose r$.
A set system $\cH$ over $X$ is \textit{$r$-uniform} if
 $\cH\ssq {X\choose r}$.
We shall use the term \emph{hypergraph} for the pair
 $(X,\cH)$ --- where $X$ is the set of \textit{vertices} and
 $\cH$ is the set of \textit{edges} or \textit{hyperedges} ---
 and also for the set system $\cH$ itself, when $X$ is understood.
The number of vertices is called the \textit{order} of $\cH$,
 and will usually be denoted by $n$.

Given a $k$-partition $\cP=(X_1,\dots,X_k)$ of $X$, we say that an $r$-set
$H\subseteq X$ \textit{crosses}  $\cP$ if $H$ intersects $\min (r,k)$
partition classes. If $r\geq k$, this means that all classes $X_i$
are intersected, whilst for $r\leq k$ the elements of the crossing
set $H$ belong to mutually distinct classes. A hypergraph $\cH$
 is said to cross $\cP$ if so does at least one of its edges $H\in \cH$.

It is a very natural problem to ask for the minimum
 number $f(n,k,r)$ of $r$-subsets (minimum number of edges in an
 $r$-uniform hypergraph), by which every
 $k$-partition of the $n$-element set $X$ is crossed.
The importance of this question is demonstrated by the fact that
 its variants have been raised by several authors independently
  in different contexts under various names:
 Sterboul in 1973 \cite{S73} (\textit{cochromatic number}, also discussed by Berge
  \cite[pp.\ 151--152]{B89}, Arocha et al.\ in 1992 \cite{ABN}
  (\textit{heterochromatic number}), and Voloshin in 1995 \cite[p.\ 43,
   Open problem 11]{V2} (\textit{upper chromatic number}, also recalled in the
    monograph \cite[Chapter 2.6, p.\ 43, Problem 2]{V}.
What is more, the formula
 $$
   f(n,2,2) = n-1
 $$
 is equivalent to the basic fact
 that every connected graph has at least $n-1$ edges
 and that this bound is tight for all $n\ge 2$.

\paragraph{Further terminology and notation.}

For a family $\frF$ of $r$-uniform hypergraphs (or graphs if $r=2$),
 and for any natural number $n$, we denote by $\tur(n,\frF)$ the
 corresponding \textit{Tur\'an number}; that is, the maximum number
 of edges in an $r$-uniform hypergraph of order $n$ that does not contain
 any subhypergraph isomorphic to any $\cF\in\frF$.
If $\frF$ consists of just one hypergraph $\cF$, we simply write
 $\tur(n,\cF)$ instead of $\tur(n,\{\cF\})$.

An $r$-uniform hypergraph $(X,\cH)$ is \textit{$r$-partite} if
 it admits a vertex partition $X_1\cup\cdots\cup X_r=X$ such that
 $|H\cap X_i|=1$ for all $H\in\cH$ and all $1\le i\le r$.
If $\cH$ consists of \textit{all} $r$-sets meeting each $X_i$ in
 precisely one vertex, then we call it a \textit{complete $r$-partite hypergraph}.

\paragraph{Earlier results.}

One can observe that a hypergraph crosses all 2-partitions of its vertex set
 if and only if it is connected.
For this reason, beyond the equation $f(n,2,2)=n-1$ mentioned above,
 we obtain that
  $$f(n,2,r)=\left\lceil\frac{n-1}{r-1}\right\rceil$$
 because this is the minimum number of edges\footnote{It is
  well known  that if $(X, \cH)$ is a
\emph{connected} hypergraph, then $\sum_{H\in \cH} \,(|H|-1) \ge
|X|-1$. The earliest source of this inequality that we have
  been able to find is Berge's classic book \cite{B},
 where Proposition 4 on page 392 is stated more generally for a
 given number of connected components.} in a
 connected $r$-uniform hypergraph of order $n$.

Let us observe further that the case of $r=2$ simply means
 graphs with at most $k-1$
 connected components, therefore
  $$f(n,k,2)=n-k+1.$$
This strong relationship with connected components, however,
 does not extend to $r>2$.

As far as we know, for $k\ge 3$ and $r\ge 3$ only the
 `diagonal case' $k=r$ of $f(n,k,r)$ has been studied
 up to now.
Below we quote the known results, using the simplified
 notation $f(n,k)$ for $f(n,k,k)$.

\begin{itemize}
\item
  $f(n,k) \ge \frac{2}{n-k+2}\,{n\choose k}$,    for every $n\ge k\ge
  3$\,\,   (\cite{S75}; later proved independently in \cite{ABN}, and
  also rediscovered in \cite{DM2000}).
\item
   $f(n,3)=\lceil \frac{n(n-2)}{3} \rceil$,   for every $n\ge 3$\,\,
   (\cite{kinai+R}; proved independently in a series of papers whose
   completing item is \cite{AT}; see also \cite{S76} for partial results).
 \item
  $f(n,n-2)={n\choose 2}-\mbox{\rm ex}(n,\{C_3,C_4\})$ \
 holds\footnote{It was quoted with a misprint in the paper \cite{BT-short}.}
  for every $n\ge 4$, where the last term is the Tur\'an number for graphs
 of girth 5 (\cite{S75}).

\end{itemize}

Although the exact value of $f(n,k)$ is not known for any $k>3$,
 its asymptotic behavior has been determined for quite a wide range of $k$.

\begin{Theorem}[\cite{BT-short}]   \label{main-uj}
 Assume\/ $n>k>2$.
\begin{itemize}
\item[$(i)$] \
$
  f(n,k) \le \frac{2}{n-1}\,{n-1\choose k} +
   \frac{n-1}{k-1}\left({n-2\choose k-2}-{n-k-1\choose k-2}\right)
$ \
  for all\/ $n$ and\/ $k$.
\item[$(ii)$] \
  $f(n,k) = (1+o(1))\, \frac{2}{k}\,{n-2\choose k-1}$ \
  for all\/ $k=o(n^{1/3})$ as\/ $n\to\infty$.
\end{itemize}
\end{Theorem}

\paragraph{Structure of the paper.} In Section~\ref{sec:2}, we first prove several preliminary results, also including an inequality for non-uniform partition-crossing hypergraphs in terms of the edge sizes. Then, we turn to uniform set systems and study the function $f(n,k,r)$ separately under the conditions $k \le r$ and $r \le k$. We prove general lower and upper bounds for $f(n,k,r)$ in both cases. In Section~\ref{sec:3}, we assume that $n-k$ and $ n-r$ are fixed while $n\rightarrow \infty$, and give asymptotically tight estimates for $f(n,k,r)$. It is worth noting that the latter problem can be reduced to Tur\' an-type problems if $k \le r$, while the same question leads us to the theory of balanced incomplete block designs if $r \le k$ is assumed.

\section{General estimates}
\label{sec:2}

Most of this section deals with uniform hypergraphs; but we shall also put
 comments on non-uniform ones which cross either all partitions or at least some
 large families of partitions.
Nevertheless the uniform systems play a central role in partition crossing,
 what will turn out already in the next subsection.

\subsection{Monotonicity}

\bp    \label{Cor-k-k'}
 For every three  integers\/ $r,k,k'$, if either \tmz
\item[$(i)$] $2\leq r\leq k\leq k'\leq n$, or
\item[$(ii)$] $2\leq k'\leq  k\leq r\leq n$ \etmz
  holds, and
  an\/ $r$-uniform hypergraph\/ $\cH$ crosses all\/ $k$-partitions of the vertex set, then\/
$\cH$ crosses all\/ $k'$-partitions, as well.
 As a consequence, for every four integers\/ $n,k,k',r$ satisfying\/ $(i)$ or\/ $(ii)$ we have
$$f(n,k,r) \geq f(n,k',r).$$ 
 \ep

\pf Assume that an $r$-uniform hypergraph $\cH$  crosses all
$k$-partitions of the vertex set $X$. Consider a
$k'$-partition $\cP'= (X_1,\dots ,X_{k'})$ of $X$.

 $(i)$ \enskip If $r\leq k\leq k'$, take the union of the last $k'-k+1$ partition classes of
$\cP'$. Due to our assumption, $\cH$ crosses
  the $k$-partition
  $\cP=(X_1, \dots ,X_{k-1}, \bigcup_{i=k}^{k'} X_i)$ obtained. Since $r\leq k$,
this means that there exists an $H\in \cH$ which contains at most
one element from
 each   partition class of $\cP$. Hence, the same $H$ and consequently, $\cH$
as well, crosses the $k'$-partition $\cP'$.

 $(ii)$ \enskip Next, assume that $k'\leq  k\leq r$ holds.
  Since the statement clearly holds for $k'=k$, we may suppose $k' < k \le n$. Then, some of the
$k'$ partition classes of $\cP'$ can be split into nonempty parts
such that a $k$-partition $\cP$ is obtained. By assumption,
 some $H\in \cH$ crosses $\cP$.
 This means that the $r$-element $H$ contains at least one element from each partition class.
 By the construction of $\cP$, $H$ contains at least one element from
  every partition class of $\cP'$;
that is, $\cH$ crosses $\cP'$.

Since the above arguments are valid for any $k'$-partition $\cP'$,
the statements follow.
 \qed

\bsk

The analogous property is valid for the other parameter of $f(n,k,r)$ as well.

\bp
 For every four  integers\/ $n,k,r,r'$, if
\tmz
\item[$(i)$] $2\leq r'\leq r\leq k\leq n$, or
\item[$(ii)$] $2\leq k \leq  r\leq r' \leq n$ holds, then \etmz
$$f(n,k,r) \geq f(n,k,r').$$
\ep
 \pf
 Consider an $r$-uniform hypegraph $(X,\cH$) of size $f(n,k,r)$ which crosses
 all $k$-partitions of the $n$-element vertex set $X$.

 $(i)$ \enskip If $r'\leq r\leq k$, then for each $H\in \cH$ choose an $r'$-element
subset $H'$ and define the $r'$-uniform set system $\cH'=\{H' \mid
H\in \cH\}$.
 Since for every $k$-partition $\cP$ there exists an $H\in \cH$ which contains
at most one element from each partition class, the same is true for
the corresponding
 $H'\in \cH'$.
Hence, $\cH'$ crosses all $k$-partitions and has at most $f(n,k,r)$
elements. This proves that $f(n,k,r) \geq f(n,k,r')$.

$(ii)$ \enskip In the other case we have $k \leq  r\leq r'$. Let
each $H\in \cH$ be extended to an arbitrary $r'$-element $H'$.  We
observe that the $r'$-uniform set system  $\cH'=\{H' \mid H\in
\cH\}$ has at most $f(n,k,r)$ elements and crosses all
$k$-partitions. Indeed, for every $k$-partition $\cP$, there exists
some $H\in \cH$ intersecting each partition class of $\cP$, and
hence the same is true for the corresponding $H'\in \cH'$. This
yields again that $f(n,k,r) \geq f(n,k,r')$ is valid. \qed

\bsk

The following corollaries show the central role of the `symmetric'
case $k=r$:

\bc   \label{c:r-to-all}
 If an $r$-uniform hypergraph $\cH$ crosses all $r$-partitions of
the vertex set $X$, then $\cH$ crosses all partitions of $X$.
\ec

Numerically, we have obtained that the function $f_{n,r}(x)=f(n,x,r)$ (where $x$ is an
integer in the range $2\leq x \leq n$) has its maximum value when
 $x=r$; and the situation is similar if $n$ and $k$ are fixed
and   $r$ is   variable; that is, the function $f_{n,k}(x)=f(n,k,x)$ attains its
maximum at $x=k$.

\bc
 \label{Cor-k-k}
For every three integers\/ $n\geq k,r \geq 2$,
    $$f(n,k,r)\leq f(n,k,k).$$ \ec

\bc
 \label{Cor-r-r}
For every three integers\/ $n\geq k,r \geq 2$,
    $$f(n,k,r)\leq f(n,r,r).$$ \ec

\subsection{Lower bound for non-uniform systems}

For hypergraphs without very small edges, we prove the following general inequality.

 \thm \label{lower-gen}
Let $k\geq 2$ be an integer, and let $(X,\cH)$ be a hypergraph of order $n$, which contains
 no edge $H\in \cH$ of cardinality smaller than $k$.
If $\cH$ crosses all $k$-partitions of $X$, then
$$
  \sum\limits_{H\in \cH} {|H| \choose k}\frac{1}{|H|-k+2}  \geq {n \choose k}\frac{1}{n-k+2}.
$$
\ethm

\pf Since $|H|\geq k$ holds for every $H\in \cH$, a $k$-partition $\cP$ of
$X$ is  crossed by $\cH$ if, and only if, there exists an edge in
$\cH$ which intersects all the $k$ partition classes of $\cP$. For every
$(k-2)$-element subset  $Y=\{x_1,\dots,x_{k-2}\}$ of $X$, define
$$\cH_Y^- = \{A \mid A\subseteq (X\setminus Y) \enskip \wedge \enskip (A\cup Y) \in \cH\}.$$
We claim that $\cH_Y^-$ is connected on $X\smin Y$.
Assume for a contradiction that it is not, and denote one of its components by $Z$.
 Consider the
$k$-partition $$\{x_1\}, \ \dots, \ \{x_{k-2}\}, \ Z, \ X\smin (Y\cup Z)$$
This is not  crossed by $\cH$ since
  a crossing set $H$ would contain all of $x_1, \dots, x_{k-2}$, moreover at least
one element from each of the last two partition classes,
 what contradicts to our assumption on disconnectivity.

Therefore,  $\cH_Y^-$ must be connected on the $(n-k+2)$-element
$X\smin Y$, and hence
$$ \sum\limits_{A \in \cH_Y^-} (|A|-1) \geq (n-k+2)-1.
$$
The corresponding inequality   holds for every $Y\in {X \choose
k-2}$. Moreover, for each edge $H\in \cH$, every $(|H|-k+2)$-element
subset of $H$ is counted in exactly one of these ${n \choose k-2}$
inequalities. Hence, we have
$$ \sum\limits_{H\in \cH} {|H| \choose k-2}(|H|-k+1) \geq {n \choose k-2}(n-k+1),
$$
which is equivalent to the assertion.
 \qed

\bsk

Beside the rather trivial hypergraph with vertex set $X$ and edge set $\cH=\{X\}$,
 which crosses every partition of $X$, the following construction
 also shows that Theorem \ref{lower-gen} is tight.

\bxm
 Let $n=|X|=2m$ be even.
 Let the edge set of $\cH$ consist of one $m$-subset $H$ of $X$
  together with $m$ mutually disjoint 2-element sets, each of which has
  precisely one vertex in $H$ and one in $X\setminus H$.
This hypergraph crosses all partitions of $X$.
Indeed, if none of the $m$ selected 2-sets crosses a partition $\cP$, then each
 class of $\cP$ meets $H$.
For this $\cH$, both sides of the inequality in Theorem~\ref{lower-gen}
  equal $\frac{n-1}{2}$ for $k=2$. (We necessarily have $k=2$,
  due to the conditions in the theorem.)
\exm

\subsection{Estimates for $k \leq r$}

The following lower bound follows immediately from Theorem~\ref{lower-gen}.

 \bc \label{lower-unif}
For every three integers\/ $n\geq r\geq k\geq 2$ the inequality
$$ f(n,k,r) \geq \frac{{n \choose k}}{{r \choose k}}\cdot \frac{r-k+2}{n-k+2}$$
holds.
  \ec

Next, we prove a general asymptotic upper bound.

\bp \label{r>k-upper}
For every two fixed integers\/ $r\geq k\geq 2$ the inequality
$$ f(n,k,r) \leq \frac{{n \choose k}}{{r \choose k}}\cdot \frac{r}{n}+ o(n^{k-1})$$
holds as\/ $n\to\infty$.
\ep

 \pf
If $k=r$, then the inequality holds also without the error term, and
 as a matter of fact, an even better upper bound on $f(n,r,r)$ is guaranteed by
 Theorem~\ref{main-uj}$(i)$.
Hence, we may suppose $r>k$.

 Consider an $n$-element vertex set $X=X'\cup \{z\}$ and an $(r-1)$-uniform
hypergraph $\cH'$ over $X'$ such that every
$(k-1)$-subset of $X'$ is covered by at least one $H'\in
\cH'$. By R\"odl's theorem \cite{Rodl}, such hypergraphs $\cH'$ of size
$$
 |\cH'|= \frac{{n-1 \choose k-1}}{{r-1 \choose k-1}} + o(n^{k-1})
$$
 exist as $n\to \infty$.

Consider now the $r$-uniform hypergraph
$$\cH=\{H'\cup \{z\} \mid H'\in \cH'\}.$$
For every $k$-partition $\cP$ we can choose a $k$-element crossing
set $A$ with $z\in A$, by picking any vertex from each of those classes of $\cP$
which do not contain $z$.
 Since $A\smin \{z\}\sst H'$ for some
$H'\in \cH'$, it follows that  $\cH$ crosses $\cP$. \qed

\bsk

We note that, beyond tight asymptotics, the
 above construction can be applied also to derive
 exact results for some restricted combinations of the parameters.

Next, we establish recursive relations to get lower bounds on $f(n,k,r)$.
Although they do not improve earlier bounds automatically, such inequalities may
 raise the possibility to propagate better estimates for larger values of the parameters
 when they are available for smaller ones.

\bp \label{rec}
If\/ $n\geq r \geq r' \geq k \geq 2$, then
$$ f(n,k,r) \geq \frac{f(n,k,r')}{f(r,k,r')}.$$
\ep
\pf Given an $n$-element vertex set $X$, consider an $r$-uniform hypergraph $\cH$
of size $f(n,k,r)$ which crosses all $k$-partitions.
Then, for each $H_j\in \cH$ construct an $r'$-uniform hypergraph $\cH_j'$
crossing all $k$-partitions of the set $H_j$. This can be done such that
$|\cH_j'|=f(r,k,r')$, hence the $r'$-uniform $\cR=\bigcup_{j=1}^{f(n,k,r)}\cH_j'$
contains at most $f(n,k,r)\cdot f(r,k,r')$ sets.

For every partition $\cP=(X_1,\dots, X_k)$, there exists some
$H_j\in \cH$ with $|X_i \cap H_j|\geq 1$ for every $1\leq i\leq k$.
Moreover, for this $j$, the system $\cH_j'$ crosses also the
$k$-partition $X_1\cap H_j, \dots, X_k\cap H_j$. Consequently, there
exists an $R\in \cH_j'\subseteq \cR$  which  intersects every class
of $\cP$. Thus, $\cR$ crosses all $k$-partitions of $X$, therefore
$$f(n,k,r)\cdot f(r,k,r')\geq f(n,k,r')$$
holds and the theorem follows. \qed

\bsk

Particularly, if $r'$ is chosen to be equal to $k$, we obtain that
$$f(n,k,r)\geq \frac{f(n,k)}{f(r,k)}.$$
Since $f(k+1,k)=k$, then
$$f(n,k,k+1)\geq \frac{f(n,k)}{k}.$$
More generally, applying Proposition~\ref{rec} repeatedly, and using the fact
 $f(i,k,i-1)=k$ that is valid for all $i>k$
  (cf.\ Proposition \ref{p:r=n-1} below), we obtain the following lower bound.
	\bc If $n\ge r \ge k\ge 2$, then 
$$f(n,k,r)\geq \frac{f(n,k)}{\prod_{i=k+1}^{r} f(i,k,i-1)}=\frac{f(n,k)}{k^{r-k}}.$$
\ec

\subsection{Estimates for $k \geq r$}

\bp \label{lower-k>r}
For every three integers $n\geq k\geq r\geq 2$ the inequality
$$ f(n,k,r) \geq \frac{{n \choose r-1}}{{k-2 \choose r-2}}\cdot \frac{n-k+2}{r(n-r+2)}$$
holds.
  \ep

\pf Consider an $r$-uniform hypergraph $\cH$ on the $n$-element
vetex set $X$, such that $\cH$ crosses all $k$-partitions.
We claim that every $(k-1)$-subset of $X$ shares at least $r-1$
 vertices with some $H\in \cH$.
Suppose for a contradiction that a set $A\in {X \choose k-1}$  intersects no $H\in \cH$
in more than $r-2$ elements.
  Then every $H\in \cH$ has at least two vertices in $X\smin A$.
Now, consider the $k$-partition whose first partition class is $X\smin A$ and the
others are singletons. This partition is not crossed by $\cH$, which is a contradiction.

Consequently, every $(k-1)$-element subset of $X$ must contain an
$(r-1)$-element subset of some $H\in \cH$. Hence,  for the `shadow'
system
$$\partial_{r-1} =\left\{B \mid  \exists H\in \cH \enskip \mbox{~s.t.} \enskip B\in{H \choose r-1}\right\},$$
the independence number must be smaller than $k-1$. Taking into
consideration the lower bound on the complementary Tur\'an number
 $T(n,k-1, r-1)={n\choose r-1} - \ext \, ( n, \cK_{k-1}^{(r-1)} )$ of complete uniform hypergraphs,
 as proved in \cite{Caen}, 
$$r\cdot |\cH|\geq |\partial_{r-1}|\geq  T(n,k-1, r-1) \geq
\frac{{n \choose r-1}}{{k-2 \choose r-2}}\cdot \frac{n-k+2}{n-r+2}$$
is obtained, from which the statement follows. \qed

\bsk

For $k$ and $r$ fixed, the lower bound gives the right order $O(n^{r-1})$, as shown by
 the following construction.

 \thm
Let $k \ge 3$, and assume that $k-2$ is divisible by  $r-2$. If $n\to\infty$, then
    $$f(n,k,r)\le \frac{2(r-2)^{r-2}}{r(k-2)^{r-2}} \, {n \choose r-1} +o(n^{r-1}) .$$
 \ethm

 \pf
Let $|X|=n$, denote $q=(k-2)/(r-2)$, and write $n'=\lceil (n-1)/q \rceil$+1.
We fix a special element $z\in X$, and partition the remaining $(n-1)$-element set $X\smin \{z\}$
 into $q$ nearly equal parts, the largest one having $n'-1$ vertices:
 $$
 X=Y_1\cup \cdots \cup Y_q \cup\{z\},   \qquad    |Y_i |=\left\lfloor
 \frac{n+i-2}{q}\right\rfloor
   \qquad \mbox{for all } \; 1\le i \le q.$$
For every set $Y_i\cup \{z\}$ we take an optimal $r$-uniform hypergraph $\cH_i$ crossing all
   $r$-partitions.
By Theorem \ref{main-uj}, we have
 $$
   |\cH_i| \le f(n',r) \le (1+o(1))\, \frac{2}{r}\,{n'-2\choose r-1} .
 $$
Here $n'-2< n/q = \frac{r-2}{k-2}\, n$, hence the binomial coefficient on the
 right-hand side is smaller than $\left( \frac{r-2}{k-2} \right)^{r-1} \! {n\choose r-1}$
Let $\cH = \cH_1 \cup \cdots \cup \cH_q$. By the estimates above, we have
 $$
   |\cH| \le  \frac{2(r-2)^{r-2}}{r(k-2)^{r-2}} \, {n \choose r-1} +o(n^{r-1})
 $$
  as $n\to\infty$.
To complete the proof, it suffices to show that $\cH$ crosses all $k$-partitions of $X$.

Let $\cP$ be any partition into $k=1+q(r-2)+1$ classes.
 One of the classes contains $z$.
By the pigeonhole principle, there is an index $i$ ($1\le i\le q$) such that, among
 the other $k-1$ classes of $\cP$ there exist at least $r-1$ which have at least one vertex in $Y_i$.
Hence we have a partition $\cP_i$ induced on $Y_i\cup\{z\}$, with some number $r'$ of classes,
 where $r'\ge r$.
Since the $r$-uniform $\cH_i$ crosses all $r$-partitions of $Y_i\cup\{z\}$,
 Corollary \ref{c:r-to-all} implies that $\cH_i$ crosses $\cP_i$, too.
That is, an $r$-set $H_i\in\cH_i$ has all its vertices in mutually distinct
 classes of $\cP_i$, which are then in distinct classes of $\cP$ as well.
Thus, $\cH$ crosses $\cP$.
 \qed

\bsk

The idea behind the construction of the above proof also yields the
 following additive upper bound.

\bp
 Suppose that all the following conditions hold:
  \begin{itemize}
    \item $n\geq k \geq r$,
    \item $n\leq 1-p+\sum_{i=1}^p n_i$,
    \item $k\leq 2-2p+\sum_{i=1}^p k_i$,
    \item $n_i\geq k_i \geq r$ for every\/  $1\leq i \leq p$.
  \end{itemize}
Then
$$f(n,k,r) \leq \sum\limits_{i=1}^p f(n_i,k_i,r).$$  \ep

\section{Asymptotics for large $k$ and $r$} 
\label{sec:3}

In this section we prove asymptotically tight estimates for $f(n,k,r)$,
 under the assumptions that the differences $s=n-k$ and $t=n-r$ are fixed and
$n\to\infty$. For this purpose, we need to consider two types of complementation ---
 one from the viewpoint of set theory, the other one analogously to graph theory.

\tmz
\item Given a hypergraph $(X, \cH)$, let $(X, \cH^c)$ denote the hypergraph of the
complements of the edges. That is, $\cH^c=\{X\smin H \mid H\in \cH\}$.
\item Given an \textit{$r$-uniform} hypergraph $(X, \cH)$, its complement $\overline{\cH}$
contains all the $r$-element subsets of $X$ which are missing from $\cH$. Formally,
$\overline{\cH}={X \choose r}\smin \cH$.
\etmz

\thm   \label{t:s<t}
Let\/ $s$ and\/ $t$ be fixed, with\/ $s\le t$, and\/ $n\to\infty$.
Then
$$
  f(n,n-s,n-t) = (1+o(1)) \frac{{n\choose s}}{{t\choose s}}.
$$
\ethm

\pf
First we prove the lower bound
 $f(n,n-s,n-t) \ge (1-o(1)) {n\choose s}/{t\choose s}$.
Suppose for a contradiction that there
exists a constant $\epsilon>0$ and an infinite sequence of
$r$-uniform hypergraphs $(X,\cH)$ with $n$ vertices and $m$
edges, edge size $r=n-t$, such that $\cH$ crosses all
 $(n-s)$-partitions of its $n$-element
vertex set $X$, but $m\le\frac{{n\choose s}}{{t\choose
s}}-\epsilon n^s$. We consider the $t$-uniform hypergraph
$\cH^c$ whose edges are the complements of the edges of $\cH$.
Since it has $m$ edges, there are at least $\epsilon {t\choose s}
n^s\ge C n^s$ distinct $s$-tuples of $X$ not covered by the edges
of $\cH^c$. Note that $C$ can be chosen as a positive absolute
constant, valid for all possible values of $n$, once we fix the triplet
$s,t,\epsilon$. We let $\cF$ to be the collection of $s$-tuples
not contained in any of the edges of $\cH^c$. Hence $|\cF|\ge C
n^s$.

Consider now the complete $s$-partite hypergraph $\cF_s$ on $2s$
vertices, each partite set having just 2 vertices. That is, the
vertex set of $\cF_s$ is $V_1\cup \cdots \cup V_s$,
 with $|V_i|=2$ for all $1\le i\le s$, and an
$s$-element set $F$ is an edge in $\cF_s$ if and only $|F\cap V_i|
= 1$ for every $i$. It is well known that the Tur\'an
number of $\cF_s$ satisfies
$$
  \ext \left( n, \cF_s \right) = o(n^s)
$$
for any fixed $s$, as $n\to\infty$. Thus, if $n$ is chosen to be
sufficiently large, $\cF$ contains a subhypergraph $\cF'$
isomorphic to $\cF_s$.

We now consider the partition $\cP$ of $X$ into $k=n-s$ classes in
which the $s$ partite sets of $\cF'$ are 2-element classes, and
the other $n-2s$ classes are singletons. By assumption, $\cH$
crosses $\cP$. It means that there exists an edge $H\in\cH$ that
meets each of the 2-element classes in at most one vertex. Let
$x_i$ be a vertex in $V_i\setminus H$ for $i=1,\dots,s$. Then
$\{x_1,\dots,x_s\} \notin \cF$, which is a contradiction to
$\{x_1,\dots,x_s\} \in \cF_s \subset \cF$,
 hence completing the proof of the lower bound.

Next, we prove the upper bound
 $
  f(n,n-s,n-t) \leq (1+o(1)) {n\choose s}/{t\choose s}
 $.
For every $n$, consider a  $t$-uniform hypergraph
$\cH_n^0$ on the $n$-element vertex set $X$,
such that each $s$-subset of $X$ is contained  in a $t$-set
$H\in \cH_n^0$.
 By R\"odl's theorem \cite{Rodl}, if $s$ and $t$ are fixed and $n\to\infty$, then
$\cH_n^0$  can be chosen such that $|\cH_n^0|={n\choose s}/{t\choose s}+o(n^s)$.

Starting with  such a  system $\cH_n^0$, we consider the hypergraph $\cH_n=(\cH_n^0)^c$
whose edge set is
$\{X\smin H \mid H\in \cH_n^0\}$. By the complementation, for $k=n-s$ and $r=n-t$, each $k$-element
subset of $X$ contains some $r$-element set $H\in \cH_n$. Then, for any $k$-partition $\cP$
of $X$, we can pick one vertex from each partition class, and this $k$-element set
has to contain an edge $H\in \cH_n$. Hence, $\cH_n$ crosses all $k$-partitions of the vertex set,
 moreover we have
$|\cH_n|=|\cH_n^0|$. This yields the claimed upper bound on
$
  f(n,n-s,n-t)
$.  \qed

\bsk

In particular, for $s=t$ we have the following consequence.
We formulate it for $s\ge 2$, because the case of $f(n,n,n)=1$ is trivial
 and the exact formula of $f(n,n-1,n-1)=n-1$ is a particular case of
 Proposition \ref{p:r=n-1} below.

 \bc For every\/ $s\ge 2$, as\/ $n\to \infty$
 $$f(n,n-s,n-s)={n\choose s}+o(n^s).$$
 \ec

To study the other range for $f(n,n-s,n-t)$, namely $s>t$,
 first we will make a simple but useful observation.
   We say that a set $T$ is a
\textit{transversal} of a partition\footnote{In fact this is the same as a transversal
 (also called vertex cover or hitting set) of the hypergraph $(X,\{X_1,\dots,X_k\})$
  in which the classes $X_i$
  of the partition are viewed as edges. This also justifies the term
  `independent set' for the complementary notion.
 }
 $\cP=(X_1,\dots,X_k)$ if $|T\cap X_i|\geq 1$ holds for every $i$.
The complement $S=X\smin T$ of a transversal $T$ is called an \textit{independent set} for $\cP$.
 This means that
$|S\cap X_i|<|X_i|$ holds for every partition class.
Let $\cI_t(\cP)$ denote the set system containing all $t$-element
independent sets for the partition $\cP$.

\bp \label{double compl}
 Let $(X,\cH)$ be an  $r$-uniform hypergraph with $|X|=n$, and assume that $k \leq r$.
Then, $\cH$ crosses all $k$-partitions of the vertex set $X$ if and only if
  for every $k$-partition $\cP$ of $X$ we have $\cI_{n-r}(\cP) \not\subseteq \overline{\cH^c}$.  
\ep

\pf For a given $k$-partition $\cP$, $\cH$ is crossing if and only if it
contains a transversal $T$ for $\cP$; that is, if $\cH^c$
contains an $(n-r)$-element independent set for $\cP$.
This equivalently means that $\overline{\cH^c}$ does not contain all elements of
$\cI_{n-r}(\cP)$. Consequently, $\cH$ crosses all $k$-partitions if  and only if
for every $k$-partition $\cP$,    $\overline{\cH^c}$ does not contain
 $\cI_{n-r}(\cP)$ as a subsystem. \qed

\bsk

Concerning $f(n,n-s,n-t)$ the case of $t=1$ is very simple.
Certainly $s=0$ means that all partition classes are singletons, hence
 $f(n,n,r)=1$ for all values of $r\le n$, also including $r=n-1$.
The situation for smaller $k$ is different.

\bp   \label{p:r=n-1}
 For every $n > k\ge 1$, we have $f(n,k, n-1)=k$.
\ep

 \pf
For $X=\{x_1,\dots,x_n\}$ define $\cH=\{X\smin\{x_i\} \mid 1\le i\le k \}$.
Consider any $k$-partition $\cP$. It either has a class with at least two
 vertices $x_i,x_j$ in the range $1\le i<j\le k$, or a class containing
 both $x_n$ and some $x_i$ with $1\le i\le k$.
Then we can choose $X\smin\{x_i\}\in\cH$, which crosses $\cP$.
Consequently, $f(n,k, n-1)\le k$.

To see the reverse inequality $f(n,k, n-1)\ge k$, without loss of generality
 we may restrict our attention to the $(n-1)$-uniform hypergraph
 $\cH^-=\{X\smin\{x_i\} \mid 1\le i\le k-1 \}$ which represents all
 $(n-1)$-uniform ones with $k-1$ edges up to isomorphism.
Then the partition
 $$
   \{x_1\}, \ \dots, \ \{x_{k-1}\}, \ \{x_k,x_{k+1},\dots,x_n\}
 $$
  is not crossed by any $H\in\cH^-$, thus $k-1$ edges are not enough.
 \qed

\bsk

The problem becomes more complicated for $t>1$.
  First we consider the case of $r=n-2$, and then a general estimate for
$k=n-s \leq n-t=r$ will be given under the assumption that $s$ and $t$ are fixed.

\bp
 For every fixed\/ $s\ge 2$,
\tmz
  \item[$(i)$] $f(n,n-s,n-2)={n\choose 2}-\tur(n,\{K_{s+1},K_{2s}-sK_2\})$;
  \item[$(ii)$] $f(n,n-s,n-2)= { \textstyle \frac{1}{2s-2} } \, n^2 +
  o(n^2)$,\,  if\/   $n\to\infty$.
 \etmz
\ep
 \pf Consider a graph $G=(V,F)$ of order $n$, which contains neither a complete graph
 $K_{s+1}$ of order $s+1$, nor a
complete graph minus a perfect matching $K_{2s}-sK_2$ on $2s$ vertices.
By the double complementation
 we obtain the $(n-2)$-uniform hypergraph $(V,\cH)=(\overline{G})^c$
  with vertex set $V$ and edge set
$$\cH=\{V\smin e \mid e\in {V \choose 2} \enskip \wedge \enskip e \notin F\}.$$
We claim that $\cH$ crosses all $(n-s)$-partitions of $V$.

First, consider a partition $\cP=(X_1,X_2, \dots ,X_{n-s})$ with at least one
partition class $|X_i|\geq 3$. We can assume without loss of generality that $|X_1|\geq 3$.
We also consider the partition $\cP'$, obtained by
removing all but one vertex from each of $X_2, \dots ,X_{n-s}$ and putting
these vertices into $X_1$. This $\cP'$ has an $(s+1)$-element class $X_1'$
and further $n-s-1$ singleton classes.
Since the class $X_1$ in $\cP$ has more than two vertices,
 every edge of $\cH$ meets $X_1$. Hence, the hypergraph $\cH$ does not cross
$\cP$ if and only if each of its edges is disjoint from at least one of the classes $X_2, \dots ,X_{n-s}$.
But then every edge is also disjoint from at least one of the singleton classes of $\cP'$, and
so $\cH$ does not cross $\cP'$ either.

Consequently, it is sufficient to ensure that $\cH$ crosses all $(n-s)$-partitions with
classes of cardinalities $(s+1, 1, \dots,1)$ and $(2,\dots,2,1, \dots 1)$,
and this will imply that $\cP$ crosses all $(n-s)$-partitions.

 An $(n-2)$-uniform  hypergraph $\cH$   crosses every partition of type
$(s+1, 1, \dots,1)$ if, and only if, for every
$(s+1)$-element subset $S$ of $V$, there exits an
edge $H\in \cH$  with $|H\cap S|=s-1$; that is, $\cH^c$ has an edge inside $S$, and equivalently,
  $G=\overline{\cH^c}$ contains no complete subgraph $K_{s+1}$.
For  the other case, $\cH$ crosses every partition of type $(2,\dots,2,1, \dots 1)$, if and only if
for every $s$ disjoint pairs of vertices there exists an edge $H$ whose
complement $\overline{H}$ contains two vertices from different pairs.
This exactly means that $G=\overline{\cH^c}$ does not contain a subgraph
$K_{2s}-sK_2$.

Consequently, an $(n-2)$-uniform $\cH$ crosses all
$(n-s)$-partitions if and only if $G$ is
$(K_{s+1},K_{2s}-sK_2)$-free. Applying the Erd\H os--Stone Theorem \cite{ErSt},
 for $s\ge 3$ this yields
$$
  f(n,n-s,n-2)={n\choose 2}-\tur(n,\{K_{s+1},K_{2s}-sK_2\})
$$
$$
  ={n\choose 2}-(1+o(1))\cdot\tur(n,K_s)
   = { \frac{1}{2s-2} } \, n^2 + o(n^2).
$$
In fact the asymptotic formula is valid also for $s=2$ because then the exclusion
 of $K_{2,2}\cong C_4$ implies that $\tur(n,\{K_{s+1},K_{2s}-sK_2\})=o(n^2)$.
\qed

\thm \label{s>t}
Let\/ $s$ and\/ $t$ be fixed, with\/ $s> t\ge 2$, and\/ $n\to\infty$.
Then,
$$
  f(n,n-s,n-t) \le  (1-c) {n\choose t}
$$
 for some constant\/ $c=c(s,t)>0$.
\ethm

\pf
Let $\cH_t$ be the complete $t$-partite hypergraph with vertex set
 $X_1\cup\cdots\cup X_t$ such that each partite class has cardinality
 $|X_i|=\lfloor n/t\rfloor$ or $|X_i|=\lceil n/t\rceil$.
We have $|\cH_t|=(1-o(1))\left(n/t\right)^t$ as $n\to\infty$,
 hence there exists a universal constant $c=c(t)>0$ such that
 $|\cH_t|\ge c {n\choose t}$ for all $n>t$.
Let $\cH=\left(\overline{\cH_t}\right)^c$.
 Then $|\cH| \le  (1-c) {n\choose t}$.

We claim that $\cH$ crosses all $(n-s)$-partitions whenever $s>t$.
 Indeed, let $\cP$ be any $(n-s)$-partition of $X$.
Consider an $s$-set $S$ obtained by deleting precisely one vertex from
 each class of $\cP$.
Since $s>t$, this $S$ contains two vertices from the same class of $\cH_t$,
 say $x',x''\in X_i$.
Therefore we can take a $t$-subset $T\sst S$ containing both $x'$ and $x''$,
 consequently $T\notin\cH_t$.
Thus, $X\smin T\in\cH$ holds, and this $X\smin T$ meets all classes of $\cP$
 because it contains all elements of $X\smin S$.
It follows that $\cH$ crosses every $\cP$, hence
$$
  f(n,n-s,n-t) \le |\cH| \le (1-c) {n\choose t} .
$$
  \qed

\end{document}